\documentclass[a4paper,10pt]{article}

\usepackage{latexsym}
\usepackage{amsmath}
\usepackage{amsthm}
\usepackage{amssymb}
\usepackage{amscd}
\usepackage{bm}
\usepackage{graphicx}
\usepackage{wrapfig}
\usepackage{fancybox}
\usepackage{enumerate}
\usepackage{url}

\newcommand{\bbr}{{\mathbb R}}
\newcommand{\bbs}{{\mathbb S}}

\newcommand{\al}{{\alpha}}

\newcommand{\gam}{{\gamma}}

\newcommand{\Del}{{\Delta}}
\newcommand{\ep}{{\epsilon}}

\newcommand{\vph}{{\varphi}}

\newcommand{\sig}{{\sigma}}

\newcommand{\re}{{\operatorname {Re}}}

\newcommand{\R}{\bbr}

\newcommand{\ovl}[1]{{\overline{#1}}}

\newcommand{\rd}{{\partial}}
\newcommand{\nab}{{\nabla}}

\newcommand{\nor}[2]{\left\| {#1} \right\|_{#2}}
\newcommand{\inp}[2]{\langle {#1} , {#2} \rangle }

\makeatletter

  \@addtoreset{equation}{section}
\makeatother

\begin{document}

\newtheorem{dfn}{Definition}[section]
\newtheorem{thm}{Theorem}[section]
\newtheorem{lem}{Lemma}[section]
\newtheorem{prop}{Proposition}[section]
\newtheorem{rmk}{Remark}[section]



\title{On uniqueness for Schr\"odinger maps with low regularity large data 
}
\author{{\sc Ikkei Shimizu} \\ Department of Mathematics, Kyoto University}
\date{\hfill}
\maketitle

\begin{abstract}
We prove that the solutions to the initial-valued problem for the 2-dimensional Schr\"odinger maps are unique in $C_tL^\infty_x \cap L^\infty_t (\dot{H}^1_x\cap\dot{H}^2_x)$. 
For the proof, we follow McGahagan's argument with improving its technical part, combining Yudovich's argument.
\end{abstract}

\section{Introduction}

We consider the initial value problem for the Schr\"odinger map equation in the two dimensional case: 

\begin{equation}\label{a1}
\left\{
\begin{aligned}
\rd_t u &= u\times \Del u  &&\text{on } \R^2\times \R, \\
u(x,0) &= u_0(x)  &&\text{on } \R^2,
\end{aligned}
\right.
\end{equation}
where $u=u(x,t)$ is the unknown function from $\R^2\times \R$ to a sphere
\begin{equation*}
\bbs^2 = \left\{ y\in \R^3 : |y|=1  \right\} \subset \R^3,
\end{equation*}
and $\times$ denotes the vector product of vectors in $\R^3$.

In the physical context, 
(\ref{a1}) is considered as a mathematical model of the evolution of magnetization vectors in ferromagnetic materials. 
(In this context, (\ref{a1}) is usually called the Landau-Lifschitz equation. For more information, see \cite{HR} for example.) 

The equation (\ref{a1}) has the energy conservation, where the energy is given by
\begin{equation*}
\mathcal{E}(u)=\int_{\R^d} \frac{1}{2} \left|\nabla u\right|^2 \, dx.
\end{equation*}
Our main theorem is the uniqueness of solutions to (\ref{a1}) in the class of $H^2$-valued solutions. The precise statement is the following:

\begin{thm}\label{T1}
(Uniqueness) Let $u^{(0)}$, $u^{(1)}\in C(I:L^\infty(\R^2))\cap L^\infty (I:\dot{H}^1(\R^2)\cap\dot{H}^2(\R^2))$ be two solutions to (\ref{a1}), where $I\subset \R$ is a time interval including $0$. 
Suppose that $u^{(0)} |_{t=0} = u^{(1)} |_{t=0}$ on $\R^2$. 
Then, $u^{(0)} = u^{(1)}$ on $\R^2\times I$.
\end{thm}

\begin{rmk}
(a) 
Theorem \ref{T1} is an improvement over that of McGahagan \cite{M} in the range of uniqueness, whose class has been the largest ever. 
Typically, McGahagan's result implies the uniqueness in $C(I:L^\infty)\cap L^\infty (I:\dot{H}^1\cap\dot{H}^{2+\ep})$ for $\ep>0$. \par
(b) By the Sobolev embedding, our result implies the uniqueness of $L^2$-localized strong solutions. 
More precisely, the solutions to (\ref{a1}) are unique in $C(I: Q+H^2)$, where $Q$ is a fixed point in $\bbs^2$.\par
(c) Theorem \ref{T1} also implies the uniqueness proved previously by the author \cite{S}. 
Under the equivariant conditions in two dimensions, he proved the uniqueness of solutions in the class $C(I:\dot{H}^1)\cap L^\infty (I:\dot{H}^2)$ 
near the harmonic family. 
Note that the auxiliary condition $C(I:L^\infty)$ follows from the embedding $\dot{H}^1\subset L^\infty$ for equivariant maps.\par 
(d) In \cite{K}, Kato proves the uniqueness of $H^1$-valued solutions to {\rm the modified Schr\"odinger map equations} under the Coulomb gauge condition. 
It is the equation which the differentiated field of the solutions to (\ref{a1}) satisfies, which is derived in \cite{NSU}. 
The $H^1$ regularity for differentiated fields corresponds to the $H^2$ regularity for original maps. 
However, note that this result does not imply the uniqueness for the original equation, as mentioned in \cite{K} and \cite{NSU}. 
\end{rmk}

The framework of the proof of Theorem \ref{T1} is mainly based on the ideas of McGahagan \cite{M}. 
In detail, we measure the difference of two solutions in a {\it geometrical} way, and show that it should be zero if the initial data coincides. 
More precisely, we consider the parallel transport of the derivatives of solutions along geodesics, and take the difference in the same tangent spaces. (See (\ref{c1}) below.) 
The advantage of this method is that we can avoid the emergence of extra derivatives on solutions, which leads to the uniqueness for $H^3$-valued solutions. 
See Proposition \ref{P2} for the mechanism of the above phenomenon. 

The differences from McGahagan's work are two points. 
First, we fixed the interval of the parameter of geodesics to $[0,1]$, instead of using the arc length parameter which is used in McGahagan's argument. 
Indeed, in \cite{M}, the quantitative estimates for geodesics (see Lemma \ref{L1}) is actually obtained by the former parameter.
However, it is sensitive to change the parameter into the later one 
when the arc length is 0, namely, the two maps coincides with each other. 
Thus we rather perform the estimates without changing the parameter.
Then, a difficult term emerges in the estimates which does not appear in \cite{M}, and 
it is necessary to exploit the geometrical structure of it. 
See the argument for the control of $Q_2$ in Section 3. \par
Second, we apply Yudovich's argument \cite{Y} to the energy method, 
Thanks to it, we can avoid the use of the Sobolev embedding for the control of $L^\infty$, which would require a little more regularity of the solutions. 
For the related result, see \cite{H} and \cite{HO} for example. \bigskip\par

The organization of the present paper is as follows. 
In Section 2, we introduce minimal geodesics and the related estimates from \cite{M}. 
We also discuss the proof of the main theorem in the same section. 
In Section 3, we give the proof of the key estimate for our proof. 
Some technical lemmas is shown in Section 4.

Finally, we introduce some notations used in the paper. 
For a function space $Y$ of space variables, 
we write $L^\infty (I: Y) := 
\{ f: I \to Y \ :\ \nor{f}{L^\infty(I:Y)} \equiv {\rm esssup}_{t\in I} \nor{f}{Y}  <\infty
\}$, $I\subset\R$. 
When $I=[0,T]$ for $T>0$, we simply write it by $L^\infty_T Y_x$.
For abbreviation, for a function space $Y$ and a vector valued function $f$, 
we write $f\in Y$ if and only if every component of $f$ belongs to $Y$. 
 We often interpret the time variable $t$ as the variable with index $0$ (see Lemma \ref{L1} for example). 
We define the commutator of two operators by $[A,B]=AB-BA$. 
For $f,g\in (L^2(\R^2))^3$, 
we define the inner product $\inp{f}{g} = \int_{\R^2} f\cdot \ovl{g} dx$. 
In our argument, we use the notation $C$ for representing a constant, and the value varies in each situations. If we intend to indicate that such $C$ depends on some specific quantity $M$, we write it by $C(M)$ for example.

\section{Main argument}

\subsection{Geodesics}
We only consider the case $I=[0,T]$ for $T>0$. The negative direction in time follows from the similar argument. 
Let $u^{(0)}$, $u^{(1)}\in C([0,T]:L^\infty)\cap L^\infty ([0,T]:\dot{H}^1\cap\dot{H}^2)$ be two solutions with $u^{(0)}|_{t=0}=u^{(1)}|_{t=1}$. 
We may choose $T$ so that $\nor{u^{(0)} -u^{(1)} }{L^\infty_{t,x}} < \pi$. 

Then, for $(x,t)\in \R^2\times [0,T]$, we consider the minimal geodesic from $u^{(0)}(x,t)$ to $u^{(1)}(x,t)$. 
In detail, we take a map $\gam (s,x,t): [0,1]\times \R^2\times [0,1]\to \bbs^2$ such that 
\begin{itemize}
\item $\gam (0,x,t) =u^{(0)} (x,t)$, $\gam (1,x,t) = u^{(1)} (x,t)$ in $(x,t)\in \R^2\times [0,T]$.
\item 
$\rd_{ss} \gam + (\gam_s \cdot \gam_s) \gamma = 0$ 
in $[0,1]\times\R^2\times [0,T]$.
\end{itemize}

The following estimates for $\gamma$ follows from the fact that the exponential map is $C^\infty$-diffeomorphism in the cut locus.  
See [Mc] in detail. 

\begin{lem}\label{L1}
(Estimates for $\gamma$)
$\gamma$ satisfies the followings: \\
(1) $|\rd_k \gam| \ \lesssim |\rd_k u^{\max} |  \quad ( k=0,1,2) $.\\
(2) $|\rd_j \rd_k \gam | \lesssim  |\rd_j \rd_k u^{\max} |+ 
 |\rd_j  u^{\max} |  |\rd_k  u^{\max} |\quad (j,k=1,2)$. 
\end{lem}

Next, we consider the parallel transport along the geodesics. 
We define the operator $X(s,\sig): T_{\gam (\sig, x,t)} \bbs^2 \to T_{\gam (s,x,t)} \bbs^2$ s.t. 
for $\xi \in T_{\gam (\sig, x,t)} \bbs^2$, $X(s,\sig ) \xi$ is the parallel transport of 
$\xi$ along $\gam (s)$. 
$X(s,\sig)$ is, by definition, the resolution operator for the following ODE:
\begin{equation*}
D_s F \equiv \rd_{s} F + (F\cdot \gam_s)\gam = 0, \quad F:\R \to \R^3.
\end{equation*}

We frequently use the following commutator estimate for $X$ and differential operators, 
which plays a central role in the argument. 
\begin{lem}\label{L2}
(Commutator) 
For $k=0,1,2$, $s,\sig\in [0,1]$ and $F(x,t) \in T_{\gam (\sig, x,t)} \bbs^2$, we have
\begin{equation}\label{a3}
\left[ D_k ,X(s,\sig) \right] F = \int_\sig^s X(s,s_1) R(\rd_s \gam , \rd_k \gam )X(s_1,\sig) Fds_1,
\end{equation}
where $R$ is the Ricci curvature tensor, defined explicitly by
\begin{equation*}
R(\xi_1, \xi_2) \xi_3 = (\xi_3\cdot \xi_2)\xi_1 - (\xi_3\cdot \xi_1) \xi_2
\end{equation*}
for $\xi_1,\xi_2,\xi_3$ in $T_p \bbs^2$ with $p\in \bbs^2$.
\end{lem}
For the proof at classical level, see \cite{M}. 
However, we need to be careful for the regularity. 
For example, in our argument in the next section, we will use Lemma \ref{L2} 
in the case when $F$ is the third derivative of $u^{(0)}$. 
Since $u^{(0)}$ is only assumed to be in $H^2$, such $F$ cannot be defined in the classical sense. 
Therefore, we have to justify the lemma above even for some \textit{distributional} vector fields. 
In Section 4, we provide a rigorous argument for this point by proving Lemma \ref{L2} in several  distributional settings (see Lemma \ref{L5}).

\subsection{Yudovich's argument}
In order to prove $u^{(0)}=u^{(1)}$, we apply Yudovich's argument to (\ref{a1}). 
The framework of the argument at formal level is the following: 
We first derive the following type of inequality for all $p>2$:
\begin{equation*}
\frac{d}{dt} ({\rm difference}) \le C p ({\rm difference})^{1-1/p}, 
\end{equation*}
where $C$ is a constant independent of $p$. By solving the above inequality, we have
\begin{equation*}
({\rm difference}) \le (C T)^p,
\end{equation*}
where we recall that $T$ is the length of time interval. 
If we choose $T$ to be sufficiently small, by taking the limit $p\to \infty$, we can conclude that the difference is $0$. 
We will proceed the above argument with caring about the singularity and regularity by 
combining penalization for instance. 
 
In our argument, we choose the difference as the following, based on \cite{M}: 
\begin{equation}\label{c1}
\nor{u^{(0)}-u^{(1)} }{L^\infty_T L^2_x}^2 + \sum_{m=1,2} \nor{X(1,0)\rd_m u^{(0)} -\rd_m u^{(1)} }{L^\infty_T L^2_x}^2.
\end{equation}

Each term is controlled in the following way:
\begin{prop}\label{P1}
Let $q= u^{(0)} -u^{(1)}$ 
and $M = \nor{\nab u^{\max}}{L^\infty_T H^1_x}$. 
Then the following estimates hold true for all $p>2$.\\
(i) $\nor{ \rd_t \nor{q }{L^2_x}^2 }{L^\infty_T} \le C(M) \nor{q}{L^\infty_T H^1_x}^{2}$.\\
(ii) $\nor{\nabla q }{L^\infty_TL^2_x} \le 
\sum_{m=1,2} \nor{X(1,0)\rd_m u^{(0)} -\rd_m u^{(1)} }{L^\infty_T L^2_x} + 
C(M)\nor{q }{L^\infty_TL^2_x}^{1-1/p}$
\end{prop}

\begin{prop}\label{P2}
Let $V_m = X(1,0)\rd_m u^{(0)} -\rd_m u^{(1)}$ for $m=1,2$, and 
$M = \nor{\nab u^{\max}}{L^\infty_T H^1_x}$. 
Then the following estimate holds true for all $p>2$.
\begin{equation}
\nor{\inp {\rd_t V_m}{ V_m}  }{L^\infty_t } \le C(M) p 
\left( \nor{q }{L^\infty_T L^2_x}^2 + \sum_{m=1,2} \nor{V_m }{L^\infty_T L^2_x}^2 \right)^{1-1/p}.
\end{equation}
\end{prop}

The main part of our argument is the proof of Proposition \ref{P2}, and we discuss it in the next section. Here we prove Proposition \ref{P1}.\bigskip\\
{\it Proof of Proposition \ref{P1}.} 
We first note that our assumption of regularity implies $u^{(0)}-u^{(1)} \in W^{1,\infty} (I : L^2)$. See Lemma \ref{L6} in Section 4 for its proof. 
Hence for a.a. $t\in I$, we have 
\begin{equation*}
\begin{aligned}
\frac{d}{dt}  \nor{q }{ L^2_x}^2 
= 2\inp{\rd_t q  }{q}
&= 2\inp{ q\times \Del u^{(0)} }{q}
+ 2\inp {u^{(1)} \times \Del q
}{q}\\
&= -2 \sum_{j=1,2} \inp {\rd_j u^{(1)} \times \rd_j q }{q}\\
&\lesssim_M \nor{q}{L^\infty_T H^1_x}^{2},
\end{aligned}
\end{equation*}
which gives (i) in the statement. For (ii), we first recall the following inequality: 
Here we recall the following inequality: For $f\in H^1$ and $r\in [2,\infty )$, we have
\begin{equation}\label{b2}
\nor{f}{L^r} \le C \sqrt{r} \nor{f}{H^1},
\end{equation}
where $C$ is independent of $p$ and $f$. 
See for instance \cite{O} for the detailed proof. \par
We decompose
\begin{equation*}
\nor{\rd_m q}{L^2_x}  \le \nor{V_m}{L^2_x} + \nor{X(1,0)\rd_m u^{(0)} - \rd_m u^{(0)}}{L^2_x}
\end{equation*}
for $m=1,2$ and $t\in[0,1]$. The second term is controlled by
\begin{equation}\label{d1}
\begin{aligned}
\nor{\int_0^1 \rd_s X(s,0) \rd_m u^{(0)}ds}{L^2_x} 
&\le \int_0^1 \nor{ \left( X(s,0) \rd_m u^{(0)} \cdot \rd_s \gam  \right)  \gam }{L^2_x} ds\\
&\le \nor{ \rd_m u^{(0)} }{L^{4p}_x} \nor{ l }{L^{4p/(2p-1)}_x} 
\end{aligned}
\end{equation}
where $l= d_{\bbs^2} (u^{(0)}, u^{(1)})$ is the distance of $u^{(0)}$ and $u^{(1)}$ on the Riemannian manifold $\bbs^2$. 
Since we have 
\begin{equation*}
|l| \le C | q |,
\end{equation*}
(\ref{d1}) is controlled by
\begin{equation*}
C \sqrt{p} \nor{\rd_m u^{(0)}}{H^1} \nor{q}{L^{4p/(2p-1)}_x} 
\le C(M) \nor{q}{L^2_x}^{1-1/p},
\end{equation*}
where we used (\ref{b2}) and the following inequality:
\begin{equation}\label{b3}
\nor{f}{L^{4p/(2p-1)}_x} \le \nor{f}{L^2_x}^{1-1/p} \nor{f}{L^4_x}^{1/p}.
\end{equation}
Thus we obtain (ii). \hfill $\square$\bigskip\\
{\it Proof of Theorem \ref{T1}.} 
We first note that $V_m\in L^\infty (I: H^1)\cap W^{1,\infty}(I:H^{-1})$ for $m=1,2$. 
We will check this regularity in Lemma \ref{L6} in Section 4. 

Let $\ep>0$ be arbitrary number. 
Let $\psi\in C^\infty (\R)$ be a cut-off function, and define the operator
$P_k = (\psi (\cdot / 2^{k}) )^{\check{\ }} * $. 
Then, for sufficiently large $k$ uniformly in $t$, we have
\begin{equation*}
\left| \nor{q}{L^2}^2 + \sum_{m=1,2} \inp{V_m}{P_k V_m} +\ep \right| >0 .
\end{equation*}
Since $P_k V_m \in C(I:H^1)$, we have
\begin{equation*}
\begin{aligned}
&\rd_t \left( \nor{q}{L^2}^2 + \sum_{m=1,2} \inp{V_m}{P_k V_m} +\ep \right)^{1/p} \\
=\frac{1}{p} & \left( \nor{q}{L^2}^2 + \sum_{m=1,2}\inp{V_m}{P_k V_m} +\ep \right)^{1/p-1} 
\left( \rd_t \nor{ q}{L^2_x}^2 + 2\sum_{m=1,2}\re\inp{\rd_t V_m}{P_kV_m} \right),
\end{aligned}
\end{equation*}
for all $t\in(0,T)$, and thus
\begin{equation*}
\begin{aligned}
&\left. \left( \nor{q}{L^2}^2 + \sum_{m=1,2} \inp{V_m}{P_k V_m} +\ep \right)^{1/p} \right|_{t=0}^{t=T}\\
= \int_0^T \frac{1}{p} & \left( \nor{q}{L^2}^2 + \sum_{m=1,2} \inp{V_m}{P_k V_m} +\ep \right)^{1/p-1} 
\left( \rd_t \nor{ q}{L^2_x}^2 + 2\sum_{m=1,2}\re\inp{\rd_t V_m}{P_kV_m} \right) dt.
\end{aligned}
\end{equation*}
By taking the limit $k\to\infty$, it follows that
\begin{equation*}
\begin{aligned}
&\left. \left( \nor{q}{L^2}^2 + \sum_{m=1,2}\nor{V_m}{L^2}^2 +\ep \right)^{1/p} \right|_{t=T}\\
= \int_0^T \frac{1}{p} & \left( \nor{q}{L^2}^2 +\sum_{m=1,2}\nor{V_m}{L^2}^2 +\ep \right)^{1/p-1} \re \left( \rd_t \nor{q}{L^2}^2 + \sum_{m=1,2} \inp{\rd_t V_m}{V_m} \right) dt,
\end{aligned}
\end{equation*}
which can be justified by Lebesgue's dominant convergence theorem.
Set 
\begin{equation}\label{b2}
G(T) = \nor{q}{L^\infty_TL^2_x}^2 + \sum_{m=1,2}\nor{V_m}{L^\infty_TL^2}^2.
\end{equation}
By Propositions \ref{P1} and \ref{P2}, we have
\begin{equation*}
(G(T) +\ep )^{1/p} \le C(M) T .
\end{equation*}
If we choose $T$ sufficiently small, $(C(M)T)^{p}\to 0$ as $p\to \infty$, and thus 
by letting $\ep\to 0$, we have
\begin{equation*}
G(T) =0,
\end{equation*}
which concludes the main theorem.

\section{Proof of Proposition \ref{P2}}
In this section, we show Proposition \ref{P2}. 
We first define $J$ by the complex structure of $\bbs^2$, which can be explicitly written as 
$J\xi =p\times \xi$ for $\xi\in T_p \bbs^2$. 
From (\ref{a1}), $V_m$ satisfies the following equation:
\begin{equation*}
D_t V_m = J \sum_{k=1}^2  D_k ^2 V_m + \sum_{\al=1}^4 R_\al,
\end{equation*}
where
\begin{equation*}
R_1 = \left[ D_t, X\right] \rd_m u^{(0)}
\end{equation*}
\begin{equation*}
R_2 = -\sum_{k=1}^2 J \left[ D_k , X\right] D_k \rd_m u^{(0)}
\end{equation*}
\begin{equation*}
R_3 = -\sum_{k=1}^2 J D_k \left[ D_k , X\right]  \rd_m u^{(0)}
\end{equation*}
\begin{equation*}
R_4 = \sum_{k=1}^2 \left\{ R(X \rd_m u^{(0)}, X \rd_k u^{(0)})X\rd_k u^{(0)} - R(\rd_m u^{(1)}, \rd_k u^{(1)}) \rd_k u^{(1)} \right\}
\end{equation*}

It suffices to show the following estimate for $m=1,2$ and $t\in I$:
\begin{equation*}
\sum_{\al=1}^4 \left|  \inp{ R_\al }{V_m} \right| \le C(M) p G(T)^{1-1/p}, 
\end{equation*}
where $G(T)$ is as in (\ref{b2}).

When $\al = 1$, we have
\begin{equation*}
\begin{aligned}
\left|  \inp{ R_1 }{V_m} \right|
&= 
\left|  
\int_0^1 \inp{ X(1,s)  R(\rd_s \gam , \rd_t \gam)  X(s,0) \rd_m u^{(0)} }{V_m} ds 
\right| \\
&\le \int_0^1 \int_{\R^2} |\rd_s \gam| |\rd_t \gam| |\rd_m u^{(0)}| |V_m| dxds\\
&\le C \int_0^1 \int_{\R^2} l |\rd_t u^{\max}| |\rd_m u^{(0)}| |V_m| dxds\\
&\le C \nor{l}{L^{8p}_x} \nor{\rd_t u^{\max}}{L^{2}_x} 
\nor{\rd_m u^{(0)}}{L^{8p}_x} \nor{V_m}{L^{4p/(2p-1)}_x}\\
&\le C(M) p \nor{q}{L^\infty H^1_x}^{1-1/p} \nor{V_m}{L^\infty_t L^2_x}^{1-1/p},
\end{aligned}
\end{equation*}
where we used Lemmas \ref{L1} and \ref{L2}, (\ref{b2}) and (\ref{b3}). 
Similarly, when $\al=2$, we have
\begin{equation*}
\begin{aligned}
\left|  \inp{ R_2 }{V_m} \right|
&= 
\left|  \sum_{k=1}^2
\int_0^1 \inp{ X(1,s)  R(\rd_s \gam , \rd_k \gam)  X(s,0) D_k \rd_m u^{(0)} }{JV_m} ds 
\right| \\
&\le C \int_0^1 \int_{\R^2} l |\nab u^{\max}| |\nab^2 u^{(0)}| |V_m| dxds\\
&\le C \nor{l}{L^{8p}_x} \nor{\nab u^{\max}}{L^{8p}_x} 
\nor{\nab^2 u^{(0)}}{L^{2}_x} \nor{V_m}{L^{4p/(2p-1)}_x}\\
&\le C(M) p \nor{q}{L^\infty H^1_x}^{1-1/p} \nor{V_m}{L^\infty_t L^2_x}^{1-1/p}.
\end{aligned}
\end{equation*}
$R_4$ is estimated in the following manner: 
\begin{equation*}
\begin{aligned}
\left|  \inp{ R_4 }{V_m} \right|
&\le C \int_{\R^2} |\nab u^{\max}|^2 (||V_1| +|V_2| )| |V_m| dxds\\
&\le C \nor{\nab u^{\max} }{L^\infty_TL^{2r}_x}^2 \nor{|V_1| + |V_2|}{L^\infty_t L^{4p/(2p-1)}_x}^2\\
&\le C(M) p \nor{|V_1| + |V_2|}{L^\infty_t L^{2}_x}^{2-2/p}.
\end{aligned}
\end{equation*}

The most difficult part is the estimate for $R_3$. 
We first consider the decomposition as follows:
\begin{equation*}
\begin{aligned}
\inp{R_3}{V_m} 
&=\sum_{k=1}^2 \int_0^1 \inp{ \left[D_k, X(1,s)\right] R(\rd_s \gam ,\rd_k \gam) X(s,0) \rd_m u^{(0)} }{JV_m} ds\\
&\hspace{40pt}+ \sum_{k=1}^2 \int_0^1 \inp{ R(D_k\rd_s \gam ,\rd_k \gam) X(s,0) \rd_m u^{(0)} }{JX(s,1))V_m} ds\\
&\hspace{40pt}+ \sum_{k=1}^2 \int_0^1 \inp{ R(\rd_s \gam , D_k\rd_k \gam) X(s,0) \rd_m u^{(0)} }{JX(s,1))V_m} ds\\
&\hspace{40pt}+ \sum_{k=1}^2 \int_0^1 \inp{ R(\rd_s \gam , \rd_k \gam) 
\left[ D_k , X(s,0) \right] \rd_m u^{(0)} }{JX(s,1))V_m} ds\\
&\hspace{40pt}+ \sum_{k=1}^2 \int_0^1 \inp{ R(\rd_s \gam , \rd_k \gam) 
X(s,0) D_k \rd_m u^{(0)} }{JX(s,1))V_m} ds\\
& \equiv \sum_{\al=1}^5 Q_\al.
\end{aligned}
\end{equation*}
Each term except $Q_2$ is estimated as follows:
\begin{equation*}
\begin{aligned}
|Q_1| &= 
\left| 
\sum_{k=1}^2 \int_0^1 \int_{s}^1 \inp{
X(1,\sig) R(\rd_s \gam , \rd_k \gam) X(\sig, s) R(\rd_s \gam, \rd_k \gam) X(s,0) \rd_m u^{(0)}
}
{JV_m} d\sig ds 
\right|\\
&\le C \int_{\R^2} l^2 |\nab u^{\max}|^3 |V_m|dx\\
&\le \nor{l}{L^\infty_T L^{8p}_x} \nor{\nab u^{\max} }{L^\infty_T L^{8p}_x} 
\nor{\nab u^{\max} }{L^\infty_T L^4_x }^2 \nor{V_m}{L^\infty_T L^{4p/(2p-1)}_x}\\
&\le C(M) p \nor{q}{L^\infty_T H^1_x}^{1-1/p} \nor{V_m}{L^\infty_t L^{2}_x}^{1-1/p}.
\end{aligned}
\end{equation*}

\begin{equation*}
\begin{aligned}
|Q_3| &\le C \int_{\R^2} l |\nab^2 u^{\max}| |\nab u^{(0)}| |V_m|\\
&\le \nor{l}{L^\infty_T L^{8p}_x} \nor{\nab^2 u^{\max} }{L^\infty_T L^{2}_x} 
\nor{\nab u^{\max} }{L^\infty_T L^{8p}_x } \nor{V_m}{L^\infty_T L^{4p/(2p-1)}_x}\\
&\le C(M) p \nor{q}{L^\infty_T H^1_x}^{1-1/p} \nor{V_m}{L^\infty_t L^{2}_x}^{1-1/p}.
\end{aligned}
\end{equation*}

\begin{equation*}
\begin{aligned}
|Q_4| &= 
\left| 
\sum_{k=1}^2 \int_0^1 \int_{s}^1 \inp{
 R(\rd_s \gam , \rd_k \gam) X(s, \sig) R(\rd_s \gam, \rd_k \gam) X(\sig,0) \rd_m u^{(0)}
}
{JV_m} d\sig ds 
\right|\\
&\le C \int_{\R^2} l^2 |\nab u^{\max}|^3 |V_m| dx\\
&\le C(M) p \nor{q}{L^\infty_T H^1_x}^{1-1/p} \nor{V_m}{L^\infty_t L^{2}_x}^{1-1/p}.
\end{aligned}
\end{equation*}

\begin{equation*}
\begin{aligned}
|Q_5| &\le C \int_{\R^2} l |\nab u^{\max}| |\nab^2 u^{(0)}| |V_m|\\
&\le C(M) p \nor{q}{L^\infty_T H^1_x}^{1-1/p} \nor{V_m}{L^\infty_t L^{2}_x}^{1-1/p}.
\end{aligned}
\end{equation*}

The argument for $Q_2$ is far involved, since we cannot control $D_k \rd_s u$ 
directly by $q$ or $V_m$. The only quantitative information of difference can be seen in the component in the 
$\rd_s \gam$ direction, that is, 
\begin{equation}\label{b5}
D_k \rd_s \gam \cdot \rd_s\gam = \frac{1}{2}\rd_k (l^2) =  l \rd_k l.
\end{equation}
It will be turned out, however, 
that the above information is sufficient to control $Q_2$ by exploiting 
the structure of Ricci tensor.\bigskip\par
Now we fix $(x,t)\in \R^2\times I$ satisfying $l(x,t)\neq 0$. 
Then we set 
\begin{equation*}
v= \frac{\rd_s \gam + J\rd_s \gam}{\sqrt{2} l},\qquad 
w= \frac{J\rd_s \gam - \rd_s \gam}{\sqrt{2} l}.
\end{equation*}
We note that $v,w$ is an orthonormal frame on $T_\gam \bbs^2$. 
By calculations, we have
\begin{equation*}
\begin{aligned}
&R(D_k \rd_s \gam, \rd_k \gam)X(s,1)\rd_m u^{(0)} \cdot JX(s,0)V_m\\
&=\left\{ 
(\rd_k \gam \cdot v) (D_k \rd_s \gam \cdot w) - (\rd_k \gam \cdot w) (D_k \rd_s \gam \cdot v)  
\right\} \\
&\hspace{20pt} \cdot 
\left\{
(X(s,0)\rd_m u^{(0)} \cdot v) (JX(s,1)V_m \cdot w) - (X(s,0)\rd_m u^{(0)} \cdot w) (JX(s,1)V_m \cdot v)
\right\}
\end{aligned}
\end{equation*}
We set
\begin{equation*}
A= (X(s,0)\rd_m u^{(0)} \cdot v) (JX(s,1)V_m \cdot w) - (X(s,0)\rd_m u^{(0)} \cdot w) (JX(s,1)V_m \cdot v).
\end{equation*}
Then $A$ is independent of $s$. Thus, we have
\begin{equation}\label{b4}
\begin{aligned}
&\int_0^1 R(D_k \rd_s \gam, \rd_k \gam)X(s,1)\rd_m u^{(0)} \cdot JX(s,0)V_m ds\\
&= A \int_0^1 \left\{ (\rd_k \gam \cdot v) (D_k \rd_s \gam \cdot w) - (\rd_k \gam \cdot w) (D_k \rd_s \gam \cdot v)  \right\}  ds \\
& = 
\frac{A}{l^2} 
\int_0^1
\left\{ (\rd_k \gam \cdot \rd_s \gam) (D_k \rd_s \gam \cdot J\rd_s \gam) - (\rd_k \gam \cdot J\rd_s \gam ) (D_k \rd_s \gam \cdot \rd_s \gam)  \right\} ds
\end{aligned}
\end{equation}
Since 
$D_k \rd_s \gam \cdot \rd_s \gam = \rd_s (D_k \gam \cdot \rd_s \gam) $, 
integration by parts gives that (\ref{b4}) is equal to
\begin{equation*}
\begin{aligned}
&\frac{2A}{l^2} \int_0^1
 (\rd_k \gam \cdot \rd_s \gam) (D_k \rd_s \gam \cdot J\rd_s \gam) ds 
- 
\left.
\frac{A}{l^2}
(\rd_k \gam \cdot J\rd_s \gam ) (D_k \gam \cdot \rd_s \gam)
\right|_{s=0}^{s=1}\\
&\equiv \frac{A}{l^2} \sum_{\al=1}^5 B_\al,
\end{aligned}
\end{equation*}
where
\begin{equation*}
B_1= 2 \int_0^1 (\rd_k \gam \cdot \rd_s \gam) (D_k \rd_s \gam \cdot J\rd_s \gam) ds ,
\end{equation*}
\begin{equation*}
B_2=  (\rd_k q \cdot \rd_s \gam|_{s=1}) (\rd_k u^{(1)} \cdot J\rd_s \gam|_{s=1})  ,
\end{equation*}
\begin{equation*}
B_3=  (\rd_k u^{(0)} \cdot \rd_s \gam|_{s=0}^{s=1}) (\rd_k u^{(1)} \cdot J\rd_s \gam|_{s=1})  ,
\end{equation*}
\begin{equation*}
B_4=  (\rd_k u^{(0)} \cdot \rd_s \gam|_{s=0}) (\rd_k q \cdot J\rd_s \gam|_{s=1})  ,
\end{equation*}
\begin{equation*}
B_5= (\rd_k u^{(0)} \cdot \rd_s \gam|_{s=0}) (\rd_k q \cdot J\rd_s \gam|_{s=0}^{s=1})  .
\end{equation*}
By using (\ref{b5}), we have
\begin{equation*}
|B_1| \le C l^2 |\rd_k l| |\rd_k u^{\max}| .
\end{equation*}
The other terms are controlled as follows:
\begin{equation*}
|B_2| +|B_4| \le C l^2 |\rd_k q| |\nab u^{\max}| 
\le C l^2 \left( |V_k| +|\rd_k u^{(0)} | l \right) |\rd_k u^{(1)}| ,
\end{equation*}
\begin{equation*}
|B_3| \le C |\rd_k u^{(0)}|  |\nab u^{(1)}| l \int_0^1 |\rd_{ss} | \gam ds
\le  C |\rd_k u^{\max}|^2 l^3 ,
\end{equation*}
\begin{equation*}
|B_5| \le C |\rd_k u^{(0)}|  |\nab u^{(1)}| l \left( \gam|_{s=0}^{s=1} \times \rd_s\gam |_{s=1} 
- \gam|_{s=0} \times \rd_s\gam|_{s=0}^{s=1}  \right) 
\le  C |\rd_k u^{\max}|^2 l^3 .
\end{equation*}
On the other hand, $A$ is bounded by
\begin{equation*}
C |\nab u^{\max}| |V_m|.
\end{equation*}
Hence, (\ref{b4}) is bounded by
\begin{equation*}
C \left\{
|u^{\max}| |\rd_k l| + |\nab u^{\max}| |V_m| + |\nab u^{\max}|^2 l
\right\}
 |\nab u^{\max}| |V_m|.
\end{equation*}
Here we recall the following estimate obtained in \cite{M}:
\begin{equation}
|\rd_k l| \le |V_k|.
\end{equation}
Writing
\begin{equation*}
\mathcal{A}(t) = \left\{ x\in \R^2 \ : \ l(x,t)=0 \right\}, 
\end{equation*}
we have
\begin{equation*}
\begin{aligned}
|Q_2| &= \left| \sum_{k=1}^2 \int_{\R^2-\mathcal{A}(t)} 
\int_0^1 R(D_k \rd_s \gam, \rd_k \gam)X(s,1)\rd_m u^{(0)} \cdot JX(s,0)V_m ds dx
\right| \\
&\le C \int_{\R^2-\mathcal{A}(t)} 
|u^{\max}|^2 (|V_1|+|V_2|)^2  dx
+C \int_{\R^2-\mathcal{A}(t)} 
|u^{\max}|^3 |V_m|^2  dx\\
&\le C(M) p \left(
\nor{V_1+|V_2|}{L^\infty_T L^2_x}^{2-2/p} 
+\nor{q}{L^\infty_T H^1_x}^{1-1/p} \nor{V_m}{L^\infty_t L^2_x}^{1-1/p}
\right),
\end{aligned}
\end{equation*}
which completes the proof. \hfill $\square$

\section{Appendix: Some technical lemmas}

We first check the following inequality from harmonic analysis:
\begin{lem}\label{L3}
For $f,g\in \mathcal{S}(\R^2)$ and $\ep>0$, we have
\begin{equation}\label{b7}
\nor{fg}{H^{-1}} \le C \nor{f}{H^{-1}} \nor{g}{L^\infty \cap\dot{H}^{1+\ep}},
\end{equation}
where $C>0$ is independent of $f,g$.
\end{lem}

\proof
By Leibniz rule (see \cite{LP} for example), we have for $\vph\in \mathcal{S} (\R^2)$, 
\begin{equation*}
\begin{aligned}
\left| \int_{\R^2} fg\vph dx \right| &\le \nor{f}{H^{-1}} \nor{g\vph}{H^{1}}\\
&\le \nor{f}{H^{-1}} \left( 
\nor{g}{L^\infty} ( \nor{\vph}{L^2} + \nor{\nab\vph}{L^2} )
+ \nor{\nab g}{L^{2/(1-\ep)}} \nor{\vph}{L^{2/\ep}}
\right)\\
&\le \nor{f}{H^{-1}} \nor{g}{L^\infty \cap\dot{H}^{1+\ep}} \nor{\vph}{H^1},
\end{aligned}
\end{equation*}
which leads to (\ref{b7}). \qedhere\bigskip\par
In particular, for $f\in H^{-1}$ and $g\in L^\infty \cap\dot{H}^{1+\ep}$ with $\ep>0$, we can define the multiplication $fg\in H^{-1}$ by the usual distributional sense,  namely, 
\begin{equation*}
\inp{fg}{\vph}_{H^{-1}, H^1} = \inp{f}{g\vph}_{H^{-1},H^1}
\end{equation*}
for $\vph\in H^1$. This is consistent with the classical definition, thus the Leibniz rule is also justified for such distributions. See \cite{S} for the related discussions. \bigskip\par
Let $\gam$ be as in the previous sections. 
We show that the parallel transport can be defined for distributional tangent vectors:

\begin{lem}\label{L4}
Let $Y$ be either 
(i) $L^\infty_T H^{-1}_x$ for $r\in (-\infty, 2)$, or (ii)  $W^{1,\infty}_TH^{-1}_x \cap L^\infty_TH^1_x$. 
Moreover, let $\sig\in[0,1]$ and $\xi\in Y^3$ satisfying $\xi \cdot \gam|_{s=\sig}=0$. 
Then a unique solution $F(s)\equiv X(s,\sig)\xi \in C^1_sY$ to the following ODE exists:
\begin{equation}\label{b8}
D_s F(s) = 0,\qquad  F(\sig ) = \xi.
\end{equation}
\end{lem}

\proof
We will show in the case $\sig=0$, where the other case can be proved similarly. 
We rewrite (\ref{b8}) in the integral form: 
\begin{equation*}
F(s) = \xi - \int_0^s \left( 
F(s')\cdot \rd_s \gam(s')
\right)
\gam (s') ds' 
\equiv \Phi (F) (s)
\end{equation*}
It suffices to show that $\Phi$ is a contraction in $C([0,s]: Y)$ for $s\ll 1$ depending only on $M$ as in Proposition \ref{P2}. In both cases of (i) and (ii), we have
\begin{equation*}
\nor{\Phi (F)}{L^\infty_s Y}
\le \nor{\xi}{Y}
+
C(M) s \nor{F}{L^\infty_s Y} 
\end{equation*}
\begin{equation*}
\nor{\Phi (F_1)- \Phi(F_2)}{L^\infty_s Y}
\le 
C(M) s \nor{F_1 -F_2}{L^\infty_s Y} 
\end{equation*}
by using Lemmas \ref{L1} and \ref{L3}, 
which proves the lemma. \qedhere

\begin{lem}\label{L5} The followings are true.\\
(i) Let $\sig\in [0,1]$ and $F\in L^\infty_TL^2_x$ with $F\cdot \gam|_{s=\sig} =0$. 
Then, (\ref{a3}) holds for $k=1,2$. \\
(ii) Let $\sig\in [0,1]$ and $F\in W^{1,\infty}_T H^{-1}_x \cap L^\infty_TH^1_x$ with $F\cdot \gam|_{s=\sig} =0$. 
Then, (\ref{a3}) holds for $k=0$.
\end{lem}

\proof 
By Lemma \ref{L4}, $\left[ D_k , X \right] F \in C^1_s L^\infty_T H^{-1}$, 
and the derivatives of the both sides of (\ref{a3}) with respect to $s$ coincides with each other. Since both sides are $0$ when $s=\sig$, the conclusion follows from the uniqueness of the solution of ODEs.\qedhere \bigskip\par
Next we check the regularity of $q$ and $V_m$.
\begin{lem}\label{L6}
The following properties hold true. \\
(i) $q\in W^{1,\infty} (I: L^2)$.\\
(ii) $V_m\in L^\infty (I: H^1)\cap W^{1,\infty}(I:H^{-1})$ for $m=1,2$.
\end{lem}

\proof
(i) From (\ref{a1}), we have $\rd_t q = q\times \nab u^{(0)} - u^{(1)} \times q \in L^\infty_T L^2_x$. In particular, it follows that $q\in W^{1,\infty}_T(L^2_x\cup L^\infty_x)\cap C_T L^\infty$. Therefore, by using Theorem 1.4.35 in \cite{CH}, we have
$
q(t, \cdot) = \int_{0}^t \rd_t q dt
$ 
for all $t\in I$, which leads to the desired conclusion. \par
(ii) By (i), it suffices to show that $X(1,0) \rd_m u^{(0)} - \rd_m u^{(0)} \in L^\infty_T H^1\cap W^{1,\infty}_T H^{-1}$. However, it immediately follows from
\begin{equation*}
X(1,0) \rd_m u^{(0)} - \rd_m u^{(0)} = - \int_0^1 (X(s,0)\rd_m u^{(0)} \cdot \rd_s \gam) \gam ds \in L^\infty_T H^1\cap W^{1,\infty}_T H^{-1},
\end{equation*}
where the last regularity property follow from Lemma \ref{L4}. \qedhere\bigskip\par
\textit{Acknowledgments} 
The author wishes to thank Yoshio Tsutsumi for the supervision, and for the many helpful comment on the present paper. 
The author also would like to thank Masayuki Hayashi for the suggestion of the use of Yudovich's argument, which leads to the improvement of the present paper. 
The author was supported by Grand-in-Aid for JSPS Fellows 18J21037.



\begin{thebibliography}{99}

	\bibitem{CH} T. Cazenave and A. Haraux, ``An introduction to semilinear evolution equations,'' 
Oxford lecture series in mathematics and its applications, 13, The Clarendon Press, New York, 1998.
	\bibitem{CSU} N. H. Chang, J. Shatah, and K. Uhlenbeck, \textit{Schr\"{o}dinger maps}, 
Comm. Pure Appl. Math., 53 (2000), 590-602.
	\bibitem{H} M. Hayashi, \textit{A note on the nonlinear Schr\"odinger equation in a general domain}, Nonlinear Anal., 173 (2018), 99-122.
	\bibitem{HO} M. Hayashi and T. Ozawa \textit{Well-posedness for a generalized derivative nonlinear Schr\"odinger equation}, J. Differ. Equations, 261 (2016), 5424-5445.
	\bibitem{HR} A. Hubert and R. Sch\"{a}fer, 
``Magnetic Domains: The Analysis of Magnetic Microstructures,'' 
Springer-Verlag, Berlin, 1998.
	\bibitem{K} J. Kato, \textit{Existence and uniqueness of the solution to 
the modified Schr\"odinger map}, Math. Res. Lett., 12 (2005), 171-186.
	\bibitem{LP} F. Linares and G. Ponce, ``Introduction to nonlinear dispersive equations. Second edition,'' Springer-Verlag, New York, 2015.
	\bibitem{M} H. McGahagan, \textit{An approximation scheme for Schr\"odinger maps}, Commun. Part. Diff. Eq., 32 (2007), 375-400.
	\bibitem{NSU} A. Nahmod, A. Stefanov, and K. Uhlenbeck, \textit{On Schr\"odinger maps}, Comm. Pure Appl. Math., 56 (2003), 114-151.
	\bibitem{O} T. Ogawa, \textit{A proof of Trudinger's inequality and its application to nonlinear Sch\"odinger equations}, Nonlinear Anal.-Theor., 14 (1990), 765-769.
	\bibitem{S} I. Shimizu, \textit{Remarks on local theory for Schr\"odinger maps near  harmonic maps}, to be appeared in Kodai Math. J., arXiv: 1806.02365.
	\bibitem{SSB} P.-L. Sulem, C. Sulem, and C. Bardos, 
\textit{On the continuous limit for a system of classical spins}, Comm. Math. Phys., 107 (1986), 431-454.  
	\bibitem{Y} V. I. Yudovich, \textit{Non-stationary flows of an ideal incompressible fluid}, \u{Z}. Vy\u{c}isl. Mat. I Mat. Fiz., 3 (1963), 1032-1066.
\end{thebibliography}
\end{document}